\documentclass{amsart}
\usepackage[all]{xypic}
\usepackage{amsmath,amsfonts}
\usepackage[pdftex]{hyperref}

\begin{document}
\bibliographystyle{plain}

\newtheorem{thm}{Theorem}[section]
\newtheorem{lem}[thm]{Lemma}
\newtheorem{prop}[thm]{Proposition}
\newtheorem{cor}[thm]{Corollary}
\newtheorem{conj}[thm]{Conjecture}
\newtheorem{mainlem}[thm]{Main Lemma}
\newtheorem{defn}[thm]{Definition}
\newtheorem{rmk}[thm]{Remark}

\def\square{\hfill${\vcenter{\vbox{\hrule height.4pt \hbox{\vrule width.4pt
height7pt \kern7pt \vrule width.4pt} \hrule height.4pt}}}$}

\newenvironment{pf}{{\it Proof:}\quad}{\square \vskip 12pt}
\newcommand{\dt}{\ensuremath{\text{det}}}
\newcommand{\U}{\ensuremath{\widetilde}}
\newcommand{\Hn}{\ensuremath{\mathbb{H}^3}}
\newcommand{\h}{{\text{hyp}}}
\newcommand{\acyl}{{\text{acyl}}}
\newcommand{\inc}{{\text{inc}}}
\newcommand{\tp}{{\text{top}}}
\newcommand{\V}{\ensuremath{\text{Vol}}}
\newcommand{\Hess}{\ensuremath{{\text{Hess} \ }}}

\title{Rigidity of minimal volume Alexandrov spaces}
\author{Peter A. Storm}

\date{April 21st, 2005}

\thanks{The author is partially supported by a National Science Foundation Postdoctoral Fellowship.}

\begin{abstract}
Let $Z$ be an Alexandrov space with curvature bounded below by $-1$ such that $Z$ is homotopy equivalent to a real hyperbolic manifold $M_\h$.  It is known that the volume of $Z$ is not smaller than the volume of $M_\h$.  If the volumes are equal, this short paper proves that the homotopy equivalence is homotopic to an isometric homeomorphism.  The main analytic tool is a theorem of Reshetnyak about quasiregular maps.
\end{abstract}
\maketitle

\section{Introduction}

The goal of this note is to prove the following version of the Besson-Courtois-Gallot theorem \cite{BCGlong} for Alexandrov spaces with curvature bounded below.

\begin{thm}
\label{Alexandrov theorem}
Let $Z$ be a compact $n$-dimensional $(n \ge 3)$ Alexandrov space with curvature bounded below by $-1$.   Let $M_\h$ be a closed oriented hyperbolic $n$-manifold.  If $f: Z \longrightarrow M_\h$ is a homotopy equivalence then
$$\text{Vol}(Z) \ge \text{Vol} (M_\h),$$
with equality if and only if $f$ is homotopic to an isometry.
\end{thm}

{\noindent}The inequality above was proven in \cite{St2}.  It remains only to prove the rigidity case.  In particular, by using the results of \cite{St2}, Theorem \ref{Alexandrov theorem} reduces to 

\begin{prop} \label{sing Reshet}
Let $N$ be a closed oriented Riemannian $n$-manifold.  Let $Z$ be a compact $n$-dimensional Alexandrov space with curvature bounded below by some $k \in \mathbb{R}$.  If $\V (Z) = \V(N)$
then a $1$-Lipschitz degree one map $Z \longrightarrow N$ is an isometric homeomorphism.
\end{prop}

{\noindent}By degree one we mean that an isomorphism is induced on $n$-dimensional integral cohomology groups.  If $Z$ is a Riemannian manifold, then the above proposition follows from a theorem of Reshetnyak (see Theorem \ref{qr maps thm} or \cite[App.C]{BCGlong}).

Theorem \ref{Alexandrov theorem} can be interpreted as evidence that Alexandrov spaces with curvature bounded below are geometrically similar to honest Riemannian manifolds.  In this vein, it may be of interest to compare their properties to those of other non-Riemannian metrics.  To facilitate this comparison, we restate Theorem \ref{Alexandrov theorem} in a slightly more general form which does not use such an explicit curvature hypothesis.

\begin{thm} \label{entropy theorem}
Let $Z$ be a compact $n$-dimensional $(n \ge 3)$ Alexandrov space with curvature bounded below by some $k \in \mathbb{R}$.   Let $M_\h$ be a closed oriented hyperbolic $n$-manifold.  Let $h(\U{Z})$ denote the volume growth entropy of the universal cover of $Z$.  If $f: Z \longrightarrow M_\h$ is a homotopy equivalence then
$$h(\U{Z})^n \, \text{Vol}(Z) \ge (n-1)^n \, \text{Vol} (M_\h),$$
with equality if and only if $f$ is homotopic to an isometry.
\end{thm}

{\noindent}Theorem \ref{entropy theorem} follows also from results of \cite{St2} and Proposition \ref{sing Reshet}.  

Stated in this form, Theorem \ref{entropy theorem} partially answers a question asked by Besson-Courtois-Gallot (see question (6) of \cite[pg.630]{BCGergodic}).  They asked if the Besson-Courtois-Gallot theorem remains true for other classes of non-Riemannian metrics, such as Finsler metrics.  Theorem \ref{entropy theorem} answers this in the affirmative for Alexandrov spaces with curvature bounded below, under the extra assumption that the source and target are connected by a homotopy equivalence, rather than simply a map of positive degree.  The Finsler case is currently still open.  There are at least two reasonable notions of volume on a Finsler manifold, and the answer may depend on which is chosen.

\vskip 6pt
This research was conducted during a visit to the Fourier Institute in Grenoble.  The author thanks his hosts for their generous hospitality.

\section{Preliminaries}

In this paper an Alexandrov space is a locally compact Alexandrov space with curvature bounded below by some $k \in \mathbb{R}$.  For more information on Alexandrov spaces, see \cite{BGP} or \cite[Ch.10]{BBI}.  Throughout we use $H^* (\cdot)$ to denote $\check{\text{C}}\text{ech}$ cohomology with coefficients in $\mathbb{Z}$.

\subsection{Quasiregular maps and degree} \label{qr maps section}

A crucial tool of this paper is the following deep result due to Reshetnyak from the theory of quasiregular maps.  


\begin{thm} \label{qr maps thm} \cite[Thm.2.16]{MV}
Let $U\subset \mathbb{R}^n$ be open.  If for some $K,c>0$, $\psi: U \longrightarrow \mathbb{R}^n$ is locally $K$-Lipschitz and $\text{Jac} \, \psi > c$ almost everywhere in $U$, then $\psi$ is an open map, there is an open subset of full measure on which $\psi$ is a local homeomorphism, and there is a $C>0$ such that for any path $\alpha \subset U$ the length of $\psi \circ \alpha$ is between $C^{-1} \cdot \text{length}(\alpha)$ and $C \cdot \text{length}(\alpha)$.
\end{thm}

\subsection{Generalized differentiable and Riemannian structures} \label{gdrs}
See \cite{OS}.  Let $X$ be a topological space, $\Omega \subseteq X$, and $n \in \mathbb{N}$.  A family $\{ (U_\phi, \phi )\}_{\phi \in \Phi}$ is called a \emph{$C^1$-atlas on $\Omega \subseteq X$} if the following hold:

(1)  For each $\phi \in \Phi$, $U_\phi$ is an open subset of $X$.

(2)  Each $\phi \in \Phi$ is a homeomorphism from $U_\phi$ into an open subset of $\mathbb{R}^n$.

(3)  $\{ U_\phi \}_{\phi \in \Phi}$ is a covering of $\Omega$.

(4)  If two maps $\phi, \psi \in \Phi$ satisfy $U_\phi \bigcap U_\psi \not= \emptyset$, then
$$\psi \circ \phi^{-1} : \phi(U_\phi \bigcap U_\psi) \longrightarrow 
        \psi (U_\phi \bigcap U_\psi)$$
is $C^1$ on $\phi(U_\phi \bigcap U_\psi \bigcap \Omega)$.
\\

A family $\{g_\phi \}_{\phi \in \Phi}$ is called a \emph{$C^{0}$-Riemannian metric} associated with a $C^1$-atlas $\{ (U_\phi, \phi) \}_{\phi \in \Phi}$ on $\Omega \subseteq X$ if the following hold:

(1)  For each $\phi \in \Phi$, $g_\phi$ is a map from $U_\phi$ to the set of positive symmetric matrices.

(2)  For each $\phi \in \Phi$, $g_\phi \circ \phi^{-1}$ is $C^0$ on $\phi(U_\phi \bigcap \Omega)$.

(3)  For any $x \in U_\phi \bigcap U_\psi, \phi,\psi \in \Phi$, we have
$$g_\psi (x) = [ d(\phi \circ \psi^{-1})(\psi(x))]^t g_\phi (x)
        [d (\phi \circ \psi^{-1}) (\psi(x))].$$

The entire reason for introducing this terminology is the following theorem.

\begin{thm}
\label{Otsu}
\cite{OS}
Let $X$ be an $n$-dimensional Alexandrov space.  Then there exists a subset $S \subset X$ of Hausdorff dimension $\le n-1$ (a set of singular points), a $C^1$-atlas on $\Omega := X \setminus S$, and a $C^0$-Riemannian metric on $\Omega$ associated with the atlas such that

(1)  The maps $\phi : U_\phi \longrightarrow \mathbb{R}^n$ of the $C^1$-atlas are locally bilipschitz.

(2)  For any $x,y \in Z$ and $\varepsilon >0$ there is a point $z \in \Omega$ such that geodesics $\alpha$ from $x$ to $z$ and $\beta$ from $z$ to $y$ have total length less than $d_Z (x,y) + \varepsilon$, and $\alpha \cup \beta \subset \Omega \cup \{ x,y\}$.  In particular, the metric structure on $\Omega$ induced from the Riemannian structure coincides with the original metric of $X$.

(3)  The Riemannian metric induces a volume element $d\text{vol}_X$ on $\Omega$.  The measure on $X$ obtained by integrating this element equals $n$-dimensional Hausdorff measure on $X$ ($S$ has zero measure).
\end{thm}

\begin{rmk}
These statements are not found in the beginning of \cite{OS}.  (1) can be found on page 651, (2) on page 654, and (3) on page 657.
\end{rmk}

{\noindent}We may therefore unambiguously define Vol$(X) := \mathcal{H}^n (X).$

We say a topological manifold is a Lipschitz manifold if there exists an atlas whose coordinate changes are bilipschitz homeomorphisms.  Condition (1) implies that an open dense subset of an Alexandrov space is a Lipschitz manifold.  In particular, condition (1) allows Rademacher's theorem to be applied to Alexandrov spaces.  This means locally Lipschitz functions on an Alexandrov space are differentiable almost everywhere.  

If $f: (U, g) \longrightarrow (V, h)$ is a continuous map between open subsets of $\mathbb{R}^n$, $g$ is a measurable Riemannian metric, $h$ is a continuous (everywhere defined) Riemannian metric, $f$ is differentiable at $p \in U_\phi$, and $g_\phi$ is defined at $p$ then we can define the unsigned Jacobian of $f$ to be
$$| \text{Jac} f | (p) :=  | \det (df)_p  | \cdot \frac{\sqrt{ \det h (f(p)) }}{\sqrt {\det g (p)}}.$$
By the last paragraph, if $f$ is Lipschitz then its unsigned Jacobian is defined almost everywhere and the change of variables theorem for Lipschitz functions can be written
$$ \int_U \chi(x) \ | \text{Jac} f |(x)  \ d \text{vol}_g (x) =  \int_V  \left( \sum_{x \in f^{-1}(y)} \chi(x) \right)  \ d \text{vol}_h (y),$$
where $\chi$ is an integrable function.  (We are implicitly using that a generic fiber will be discrete.  See \cite[Sec.3.4.3]{EG} for the change of variable theorem for Lipschitz functions.)  For later use we now prove a weak change of variables theorem for proper Lipschitz maps between Lipschitz manifolds.  The proof is straightforward, but we include it for completeness.

\begin{lem} \label{change of vars}
Let $f: M \longrightarrow N$ be a locally Lipschitz proper map between Lipschitz manifolds.  Let $g$ be an almost everywhere defined measurable Riemannian metric on $M$, and let $h$ be a continuous (everywhere defined) Riemannian metric on $N$.  Then
$$\int_M  | \text{Jac} f  | \ d\text{vol}_g  \ge  \int_N  \# \{ f^{-1} (y) \} \ d\text{vol}_h (y).$$
\end{lem}

\begin{pf}
Let $\mathcal{D} \subset N$ be the full measured set of points with discrete fiber.  For each point $q \in \mathcal{D}$ pick an open neighborhood $V_q \subset N$ of $q$ such that $V_q$ is contained in a single coordinate chart of $N$, and each component of $f^{-1} (V_q) \subset M$ is contained in a single coordinate chart of $M$.  (This is where we use the assumption that $f$ is proper.)  The collection $\{ V_q \}_{q \in \mathcal{D}}$ covers the full measured open set $\cup_{q \in \mathcal{D} } V_q$.  Pick a locally finite subcover $\{ V_i \}$ of the same set.  Pick a partition of unity $\{ \eta_j \}$ subordinate to the cover $\{ V_i \}$.  We may now apply the change of variables theorem in Euclidean space to conclude that for any $i$ and $j$
$$\int_{V_i} \eta_j (y) \  \# \{f^{-1}(y) \} \ d\text{vol}_h (y) = \int_{f^{-1}(V_i)} | \text{Jac} f | \, (\eta_j \circ f) \ d\text{vol}_g.$$
The lemma follows by summing over all $i,j$ and noting that $\cup_i f^{-1}(V_i) \subseteq M$.
\end{pf}

\section{The Proof} \label{appC thm}

As stated in the introduction, the main theorem will follow from

\vskip 6pt
\noindent\textbf{Proposition \ref{sing Reshet}.} \itshape
Let $N$ be a closed oriented Riemannian $n$-manifold.  Let $Z$ be a compact $n$-dimensional Alexandrov space with curvature bounded below by some $k \in \mathbb{R}$.  If $\V (Z) = \V(N)$
then a $1$-Lipschitz degree one map $Z \longrightarrow N$ is an isometric homeomorphism.
\normalfont \vskip 6pt

{\noindent}Recall that by degree one we mean the map induces an isomorphism $H^n (N) \rightarrow H^n(Z)$.

\emph{Proof:}
Recall that $\Omega \subseteq Z$ is a subset of full measure possessing a $\mathcal{C}^1$-atlas and a $\mathcal{C}^0$-Riemannian metric (see Section \ref{gdrs}).  Let $\mathcal{U} \subseteq Z$ be a connected open subset of full measure containing $\Omega \subseteq Z$ which is locally bilipschitz to Euclidean space.  $\mathcal{U}$ is obtained as the union of the sets $U_\phi$ of Section \ref{gdrs}, and the locally bilipschitz charts $\phi : U_\phi \longrightarrow \mathbb{R}^n$ put a Lipschitz structure on the manifold $\mathcal{U}$.  Define $\partial \mathcal{U} := \overline{\mathcal{U}} - \mathcal{U} = Z - \mathcal{U} \subset Z$.

Let $\Phi: Z \longrightarrow N$ be a $1$-Lipschitz degree one map.  

\begin{lem} \label{topological lem}
The topological manifold $\mathcal{U} \subseteq Z$ is orientable, the inclusion $(Z, \emptyset) \subseteq (Z, \partial \mathcal{U})$ induces an isomorphism 
$$\eta: H^n (Z, \partial \mathcal{U}) \longrightarrow H^n (Z) \cong H^n (N) \cong \mathbb{Z},$$
and for any open topological $n$-cell $\mathcal{O} \subseteq \mathcal{U}$ inclusion 
$$(Z , \partial \mathcal{U}) \subseteq (Z, Z- \mathcal{O})$$ induces an isomorphism 
$$\iota_{\mathcal{O}}^* : H^n (Z, Z- \mathcal{O}) \longrightarrow H^n (Z, \partial \mathcal{U}).$$
\end{lem}
{\noindent}Keep in mind that \emph{a priori} $Z$ may not be a topological manifold.

\begin{pf}
The pair $(Z, \partial \mathcal{U})$ has a long exact sequence on cohomology.  At dimension $n$ it is
$$H^n (Z, \partial \mathcal{U}) \stackrel{\eta}{\longrightarrow} H^n (Z) \longrightarrow H^n (\partial \mathcal{U}).$$
Using the fact that $\Phi$ is degree one, we know $H^n (Z) \cong \mathbb{Z}$.  Since $\partial \mathcal{U}$ is contained in the singular set $S$ of $Z$, and the Hausdorff dimension of $S$ is at most $n-1$, it follows that the topological dimension of $\partial \mathcal{U} \subseteq S$ is at most $n-1$ \cite[VII.4]{HW}.  Therefore $H^n (\partial \mathcal{U}) = 0$, implying that $\eta$ is surjective.  The pair $(Z, \partial \mathcal{U})$ is a relative topological manifold.  So either $H^n (Z, \partial \mathcal{U}) \cong \mathbb{Z}$ and $\mathcal{U}$ is orientable, or $H^n (Z, \partial \mathcal{U}) \cong \mathbb{Z}/2\mathbb{Z}$ and $\mathcal{U}$ is not orientable \cite[XI.6.8]{ES}.  The lemma now follows from the surjectivity of $\eta$ and an application of \cite[XI.6.8]{ES}.
\end{pf}

We will use open sets in $\mathcal{U}$ with particularly nice standard properties.  Let us call an open set $\mathcal{O} \subset \mathcal{U}$ a \emph{good} open set if its closure lies in a single coordinate chart $U_\phi$ of $\mathcal{U}$, and $(\overline{\mathcal{O}}, \partial \mathcal{O})$ is homeomorphic as a pair to $(\overline{ \mathbb{B}}, \partial \mathbb{B})$.  (Here $\partial \mathcal{O} := \overline{\mathcal{O}}- \mathcal{O}$ and $\mathbb{B}$ is the open unit ball in $\mathbb{R}^n$.)

Let us now see how to use Lemma \ref{topological lem} to orient good open sets in $\mathcal{U}$.  Fix for once and for all a generator $1_N \in H^n (N)$.  Then the isomorphisms
$$H^n (Z, \partial \mathcal{U}) \stackrel{\eta}{\longrightarrow} H^n (Z) \stackrel{\Phi^*}{\longleftarrow} H^n (N)$$
determine an element $1_\mathcal{U} \in H^n (Z, \partial \mathcal{U})$.  Pick a good open set $\mathcal{O} \subset \mathcal{U}$.  Then the isomorphism $\iota_{\mathcal{O}}^*$ determines an element of $H^n (Z, Z- \mathcal{O})$.  By excision
$$H^n (Z, Z- \mathcal{O}) \stackrel{\sim}{\longleftarrow} H^n (\overline{\mathcal{O}}, \partial \mathcal{O}) \cong  \mathbb{Z}.$$ 
Pull $\iota_{\mathcal{O}}^* (1_\mathcal{U})$ back by the excision isomorphism.  The resulting element $1_\mathcal{O} \in H^n (\overline{\mathcal{O}}, \partial \mathcal{O})$ determines an orientation of $\mathcal{O}$.

Keeping this notation, pick a point $p \in \mathcal{O}$ where $| \text{Jac} \Phi |$ is defined.  We can now give it a sign in the natural way: pick a basis of $T_p \mathcal{O}$ which is positively oriented with respect to $1_\mathcal{O}$, and give the Jacobian a sign based on whether or not $d\Phi$ pushes that basis forward to a positively or negatively oriented basis of $T_{\Phi(p)} N$ with respect to $1_N$.

Note that $\Phi$ is volume non-increasing and $\V (Z) = \V(N)$.  From this it follows that $\Phi$ is volume preserving.  In particular both $\Phi (\partial \mathcal{U}) \subset N$ and $\Phi^{-1}(\Phi(\partial \mathcal{U})) \subset Z$ have measure zero.  Note that the restricted map 
$$\Phi : \mathcal{U} - \Phi^{-1}(\Phi(\partial \mathcal{U})) \longrightarrow N - \Phi(\partial \mathcal{U})$$
is a proper locally Lipschitz map between open sets of full measure.

\begin{lem}
$| \text{Jac} \Phi | = 1$ almost everywhere and $\# \{ \Phi^{-1} (y) \} = 1$ for almost every $y \in N$.
\end{lem}
\begin{pf}

By the change of variables theorem (see Lemma \ref{change of vars})
$$\text{Vol}(Z) \ge  \int_{\mathcal{U}- \Phi^{-1}(\Phi(\partial \mathcal{U}))} | \text{Jac} \Phi | \ d \text{vol}_Z \ge 
			\int_{ N - \Phi(\partial \mathcal{U})} \# \{ \Phi^{-1}(y) \} \ d\text{vol}_N \ge \text{Vol}(N).$$
Recalling that $\text{Vol}(Z)= \text{Vol}(N)$, the lemma follows immediately from the above inequalities.
\end{pf}

In order to apply Reshetnyak's Theorem \ref{qr maps thm}, we must nail down the sign of the Jacobian almost everywhere.  Recall that the degree of a locally Lipschitz proper map between smooth manifolds is the oriented cardinality of a generic fiber \cite[p.383, cor.4.1.26]{Fed}.  (In this statement it is assumed that the target manifold is connected.)

\begin{lem}
$\text{Jac} \Phi = 1$ almost everywhere.
\end{lem}
\begin{pf}
Pick a good neighborhood $V_1 \subseteq N - \Phi(\partial \mathcal{U})$ containing a point $y$ such that $\Phi^{-1}(y)$ is a single point.  (Such points $y$ are generic.)  Pick a good neighborhood $\mathcal{O} \subseteq \Phi^{-1}(V_1)$ containing $x := \Phi^{-1}(y)$.  Pick a good neighborhood $y \in V_2 \subseteq V_1$ sufficiently small such that $\overline{\Phi^{-1}(V_2)} \subset \mathcal{O}$.  

Implicitly using the Lipschitz coordinate charts, consider the Lipschitz map between closed Euclidean balls given by
$$\Phi : ( \overline{\mathcal{O}}, \partial \mathcal{O}) \longrightarrow (\overline{V_1}, \overline{V_1} - V_2 ).$$
The degree of this map can be determined topologically by chasing the following commutative diagram.

$$\xymatrix{
{H^n (Z, \emptyset)}   &&  \ar[ll]_(.4){\Phi^*}  {H^n (N, \emptyset)}   \\
{H^n (Z, \partial \mathcal{U})}  \ar[u] &&   \\
{H^n (Z, Z - \mathcal{O})}  \ar[d]_{\text{excision}} \ar[u] &&  \ar[ll]_{\Phi^*} {H^n (N, N - V_2)} \ar[uu] \ar[d]_{\text{excision}}  \\
{H^n (\overline{\mathcal{O}}, \partial \mathcal{O})}   && \ar[ll]_{\Phi^*} {H^n (\overline{V_1}, \overline{V_1}- V_2)}  }$$

{\noindent}Here every arrow is an isomorphism.  Moreover, as every vertical arrow is induced by topological inclusion, the diagram obviously commutes.  Beginning with an orientation $1_N \in H^n (N, \emptyset)$, the orientation $1_\mathcal{O} \in H^n (\overline{\mathcal{O}}, \partial \mathcal{O})$ is determined by going around the edge of the diagram counterclockwise.  The orientation of $(\overline{V_1}, \overline{V_1} - V_2)$ is obtained by going down from $H^n (N, \emptyset)$.  From commutativity it follows that the degree of 
$$\Phi : ( \overline{\mathcal{O}}, \partial \mathcal{O}) \longrightarrow (\overline{V_1}, \overline{V_1} - V_2 )$$
is $1$.  Therefore the oriented cardinality of a generic fiber over $V_2$ is one.  Using the previous lemma it follows that $\text{Jac} \Phi$ is $1$ almost everywhere in $\Phi^{-1} (V_2)$.  The chosen point $y \in V_2$ was generic and $\Phi$ preserves measure.  The lemma follows from these facts.
\end{pf}

Using the fact that $\mathcal{U} \subseteq Z$ is locally bilipschitz to Euclidean space, we may repeatedly apply Theorem \ref{qr maps thm} on small neighborhoods to conclude that the map $\Phi: \mathcal{U} \longrightarrow N$ is open, and that for some open set of full measure $\mathcal{U}' \subseteq \mathcal{U}$ the restricted map $\Phi: \mathcal{U}' \longrightarrow N$ is a local homeomorphism onto its image.  

\begin{lem} \label{homeo lem}
$\Phi: \mathcal{U} \longrightarrow N$ is a homeomorphism onto its image.
\end{lem}
\begin{pf}
Suppose there exist distinct points $z_1, z_2 \in \mathcal{U}$ such that $\Phi(z_1) = \Phi (z_2) = p \in N$.  The map $\Phi$ is open.  Therefore there exists an open neighborhood $\mathcal{O}_i$ of $z_i$ such that $\mathcal{O}_1 \cap \mathcal{O}_2 = \emptyset$ and $\Phi(\mathcal{O}_1)= \Phi(\mathcal{O}_2)$ is an open neighborhood of $p$.  Since $\Phi(Z - \mathcal{U}') \subset N$ has measure zero, $\Phi(\mathcal{O}_i) - \Phi(Z-\mathcal{U}')$ is nonempty and open.  For any point $q \in \Phi(\mathcal{O}_i) - \Phi(Z - \mathcal{U}')$, the fiber $\Phi^{-1}(q)$ contains at least two points in $\mathcal{U}$.  This contradicts the fact that the generic fiber is a single point.
\end{pf}

Using Lemma \ref{homeo lem}, the conclusion of Theorem \ref{qr maps thm} automatically promotes $\Phi : \mathcal{U} \longrightarrow N$ to a locally bilipschitz map.  In fact, since $\Phi$ is $1$-Lipschitz and $\text{Jac} \Phi = 1$ almost everywhere, it follows that the map of tangent spaces $d (\Phi^{-1})_p$ is $1$-Lipschitz for almost every $p \in N$.  A standard path-deformation argument using Fubini's theorem (see, for example, the proof of \cite[Lem.7.8]{BCGlong} or \cite[Lem.8.4]{St2}) can now be used to show that $\Phi^{-1}$ is locally $1$-Lipschitz.  Therefore $\Phi|_\mathcal{U}$ is a local isometry.

It remains only to deal with points in $\partial \mathcal{U} \subset Z$.  Fix $\delta >0$ such that any $\delta$-ball in $N$ is strictly convex.  Pick $x,y \in Z$ such that $d_Z (x,y) < \delta/2$.  Let $\alpha_\varepsilon \subset \mathcal{U} \cup \{ x,y \}$ be a once broken piecewise geodesic path formed by concatenating a geodesic from $x$ to $z_\varepsilon$ and from $z_\varepsilon$ to $y$ (for some $z_\varepsilon \in \mathcal{U}$), with total length less than $d_{Z} (x,y) + \varepsilon$.  Such a path $\alpha_\varepsilon$ exists by Theorem \ref{Otsu}.  Since $Z$ is an Alexandrov space with curvature $\ge -1$, a comparison geometry argument shows that the angle formed by $\alpha_\varepsilon$ at $z_\varepsilon$ goes to $\pi$ (i.e. $\alpha_\varepsilon$ becomes more straight) as $\varepsilon \rightarrow 0$.  $\Phi|_\mathcal{U}$ is a local isometry, implying the path $\Phi \circ \alpha \subset N$ is also a piecewise geodesic path with small curvature at its one point of nondifferentiability.  $\Phi$ is $1$-Lipschitz, so for $\varepsilon < \delta/2$ the path $\Phi \circ \alpha$ is contained in a strictly convex $\delta$-ball in $N$.  In particular, by taking $\varepsilon \rightarrow 0$ we see that $d_N (\Phi(x), \Phi(y)) = d_{Z} (x,y)$.  Therefore $\Phi$ is locally an isometric homeomorphism.  Since $\Phi$ is degree one, it is a global isometry.  This completes the proof of Proposition \ref{sing Reshet}. \square \vskip 12pt


We can now easily prove the main theorem of this paper.

\vskip 6pt
\noindent\textbf{Theorem \ref{Alexandrov theorem}.} \itshape
Let $Z$ be a compact $n$-dimensional $(n \ge 3)$ Alexandrov space with curvature bounded below by $-1$.   Let $M_\h$ be a closed oriented hyperbolic $n$-manifold.  If $f: Z \longrightarrow M_\h$ is a homotopy equivalence then
$$\text{Vol}(Z) \ge \text{Vol} (M_\h),$$
with equality if and only if $f$ is homotopic to an isometry. \normalfont \vskip 6pt

\begin{pf}
The inequality of Theorem \ref{Alexandrov theorem} was proven in \cite{St2}.  So assume that $\text{Vol}(Z) = \text{Vol}(M_\h)$.  Then the natural map $F: Z \longrightarrow M_\h$ given by \cite[Prop.5.5]{St2} is Lipschitz, homotopic to $f$, and $dF$ is an infinitesimal isometry almost everywhere.  By perturbing paths in the open set $\mathcal{U} \subseteq Z$ and using Fubini's theorem, a standard argument (see \cite[Lem.8.4]{St2}) can be used to show that $F: Z \longrightarrow M_\h$ is $1$-Lipschitz.  Proposition \ref{sing Reshet} then implies that $F$ is an isometry.
\end{pf}

{\noindent}Theorem \ref{entropy theorem} is proven identically.  (See \cite[Sec.5]{St2} and note that an Alexandrov space with curvature bounded below is an ``almost every Riemannian metric space.'')

\bibliography{rmvas.biblio}

\end{document}